\documentclass[12pt]{article}%
\usepackage{amssymb}
\usepackage{graphicx}
\usepackage[intlimits]{amsmath}
\usepackage{latexsym}
\usepackage{amsfonts}%
 \makeatletter
\newfont{\footsc}{cmcsc10 at 8truept}
\newfont{\footbf}{cmbx10 at 8truept}
\newfont{\footrm}{cmr10 at 10truept}
\makeatother 
\newtheorem{theorem}{Theorem}

\newtheorem{definition}[theorem]{Definition}
\newtheorem{example}[theorem]{Example}

\newtheorem{lemma}[theorem]{Lemma}

\newenvironment{proof}[1][Proof]{\noindent{\textbf {#1}  }}  {\hfill$\Box$\bigskip}

\begin{document}

\title{An Abstract Regularity Lemma}
\author{B. Bollob\'{a}s, V. Nikiforov\\Memphis}
\maketitle

\begin{abstract}
We extend in a natural way Szemer\'{e}di's Regularity Lemma to abstract
measure spaces.

\end{abstract}

\section{Introduction}

In this note we extend Szemer\'{e}di's Regularity Lemma (SRL) to abstract
measure spaces. Our main aim is to find general conditions under which the
original proof of Szemer\'{e}di still works. Another extension of SRL to
probability spaces was proved by Tao \cite{Tao06}, but his results do not
imply our most general result, Theorem \ref{gmaint}. To illustrate that our
approach has some merit, we outline several applications. Some of these
applications seem to be tailored to our approach: in particular, we are not
aware of any alternative proofs.

Our notation follows \cite{Bol98}.

\subsection{Measure triples}

A finitely additive measure triple or, briefly, \emph{a measure triple}
$\left(  X,\mathbb{A},\mu\right)  $ consists of a set $X,$ an algebra
$\mathbb{A}\subset2^{X}$, and a complete, nonnegative, finitely additive
measure $\mu$ on $\mathbb{A}$ with $\mu\left(  X\right)  =1$. Thus,
$\mathbb{A}$ contains $X$ and is closed under finite intersections, unions and
differences; the elements of $\mathbb{A}$ are called \emph{measurable subsets}
of $X.$\bigskip

Here are some examples of measure triples.

\begin{example}
\label{exi}Let $k,n\geq1,$ write $2^{\left[  n\right]  ^{k}}$ for the power
set of $\left[  n\right]  ^{k},$ and define $\mu^{k}$ by $\mu^{k}\left(
A\right)  =\left\vert A\right\vert /n^{k}$ for every $A\subset\left[
n\right]  ^{k}$. Then $\left(  \left[  n\right]  ^{k},2^{\left[  n\right]
^{k}},\mu^{k}\right)  $ is a measure triple.
\end{example}

Note that there is a one-to-one mapping between undirected $k$-graphs on the
vertex set $\left[  n\right]  $ and subsets $G\subset2^{\left[  n\right]
^{k}}$ such that if $\left(  v_{1},\ldots,v_{k}\right)  \in G,$ then $\left\{
v_{1},\ldots,v_{k}\right\}  $ is a set of size $k$ and $G$ contains every
permutation of $\left(  v_{1},\ldots,v_{k}\right)  $. In view of this, we
shall consider subsets of $2^{\left[  n\right]  ^{k}}$ as labelled directed
$k$-graphs (with loops) on the vertex set $\left[  n\right]  .$

\begin{example}
\label{exmp}Let $k\geq2,$ and let $X_{1},\ldots,X_{k}$ be finite nonempty
disjoint sets. Write $2^{X_{1}\times\cdots\times X_{k}}$ for the power set of
$X_{1}\times\cdots\times X_{k},$ and define $\mu^{k}$ by $\mu^{k}\left(
A\right)  =\left\vert A\right\vert /\left(  \left\vert X_{1}\right\vert
\cdots\left\vert X_{k}\right\vert \right)  $ for every $A\subset X_{1}%
\times\cdots\times X_{k}$. Then $\left(  X_{1}\times\cdots\times
X_{k},2^{X_{1}\times\cdots\times X_{k}},\mu^{k}\right)  $ is a measure triple.
\end{example}

We shall consider subsets of $2^{X_{1}\times\cdots\times X_{k}}$ as labelled
$k$-partite $k$-graphs with vertex classes $X_{1},\ldots,X_{k}.$

\begin{example}
\label{exl}Let $k\geq1,$ and let $\mathbb{B}^{k}$ be the algebra of the Borel
subsets of the unit cube $\left[  0,1\right]  ^{k};$ write $\mu^{k}$\ for the
Lebesgue measure on $\mathbb{B}^{k}.$ Then $\left(  \left[  0,1\right]
^{k},\mathbb{B}^{k},\mu^{k}\right)  $ is a measure triple.
\end{example}

\subsection{SR-systems}

Let us introduce the main objects of our study, SR-systems: measure triples
with a suitably chosen semi-ring. Here SR stands for \textquotedblleft
Szemer\'{e}di Regularity\textquotedblright\ rather than \textquotedblleft
semi-ring\textquotedblright.

Recall that a set system $\mathbb{S}$ is a semi-ring if it is closed under
intersection and for all $A,B\in\mathbb{S}$, the difference $A\backslash B$ is
a disjoint union of a finite number of members of $\mathbb{S}$.

A semi-ring $\mathbb{S}$ is called $r$\emph{-built} if for all $A,B\in
\mathbb{S}$, the difference $A\backslash B$ is a disjoint union of at most $r$
members of $\mathbb{S}$; we say that $\mathbb{S}$ is \emph{boundedly built} if
it is $r$-built for some $r$.

An \emph{SR-system} is a quadruple $\left(  X,\mathbb{A},\mu,\mathbb{S}%
\right)  $, where $\left(  X,\mathbb{A},\mu\right)  $ is a measure triple and
$\mathbb{S}\subset\mathbb{A}$ is\ a boundedly built semi-ring. Clearly the
quadruple $\left(  X,\mathbb{A},\mu,\mathbb{A}\right)  $ is the simplest
example of an SR-system based on the measure triple $\left(  X,\mathbb{A}%
,\mu\right)  .$

For the rest of the section, let us fix an SR-system $\left(  X,\mathbb{A}%
,\mu,\mathbb{S}\right)  $.

Given a set system $\mathbb{X}$ and $k\geq1,$ let $\mathbb{X}^{\left\langle
k\right\rangle }$ be the collection of products of $k$ elements of
$\mathbb{X}$ any two of which are either disjoint or coincide, i.e.,
\[
\mathbb{X}^{\left\langle k\right\rangle }=\left\{  A_{1}\times\cdots\times
A_{k}:A_{i}\in\mathbb{X}\text{ and }A_{i}\cap A_{j}=\varnothing\text{ or
}A_{i}=A_{j}\text{ for all }i,j\in\left[  k\right]  \right\}  .
\]

The proof of the following lemma is given in Section \ref{sprf}.

\begin{lemma}
\label{le1}The set system $\mathbb{S}^{\left\langle k\right\rangle }$ is a
boundedly built semi-ring.\bigskip
\end{lemma}

This assertion is used in the following general construction .

\begin{example}
\label{exq copy(1)} For $k\geq1$, set
\[
\mathbb{A}^{k}=\left\{  A_{1}\times\cdots\times A_{k}:A_{i}\in\mathbb{A}\text{
for all }i\in\left[  k\right]  \right\}  .
\]
Write $\mathcal{A}({\mathbb{A}}^{k})$ for the algebra generated by the set
system $\mathbb{A}^{k}$, and $\mu^{k}$ for the product measure on
$\mathcal{A}(${$\mathbb{A}$}$^{k}).$ The quadruple $\left(  X^{k}%
,\mathcal{A}({\mathbb{A}}^{k}),\mu^{k},\mathbb{S}^{\left\langle k\right\rangle
}\right)  $ is an SR-system.
\end{example}

Let us see three particular examples of the above construction.

\begin{example}
For $k\geq1$ set $\mathcal{G}^{k}\left(  n\right)  =\left(  \left[  n\right]
^{k},2^{\left[  n\right]  ^{k}},\mu^{k},\left(  2^{\left[  n\right]  }\right)
^{\left\langle k\right\rangle }\right)  $, where $\left(  \left[  n\right]
^{k},2^{\left[  n\right]  ^{k}},\mu^{k}\right)  $ is the measure triple
defined in Example \ref{exi}, and $\left(  2^{\left[  n\right]  }\right)
^{\left\langle k\right\rangle }$ is the set of all products of $k$ subsets of
$\left[  n\right]  $ any two of which are either disjoint or coincide$.$
\end{example}

\begin{example}
For $k\geq1$ set $\mathcal{B}^{k}=\left(  \left[  0,1\right]  ^{k}%
,\mathbb{B}^{k},\lambda^{k},\mathbb{B}^{\left\langle k\right\rangle }\right)
$, where $\left(  \left[  0,1\right]  ^{k},\mathbb{B}^{k},\lambda^{k}\right)
$ is the measure triple defined in Example \ref{exl}, and $\mathbb{B}%
^{\left\langle k\right\rangle }$ is the set of all products of $k$ Borel
subsets of $\left[  0,1\right]  $ any two of which are either disjoint or coincide.
\end{example}

\begin{example}
For $k\geq1$ set $\mathcal{BI}^{k}=\left(  \left[  0,1\right]  ^{k}%
,\mathbb{B}^{k},\lambda^{k},\mathbb{I}^{\left\langle k\right\rangle }\right)
$, where $\left(  \left[  0,1\right]  ^{k},\mathbb{B}^{k},\lambda^{k}\right)
$ is the measure triple defined in Example \ref{exl}, and $\mathbb{I}%
^{\left\langle k\right\rangle }$ is the set of all products of $k$ intervals
$\left[  a,b\right)  \subset\left[  0,1\right]  $ any two of which are either
disjoint or coincide.
\end{example}

Here is another general SR-system.

\begin{example}
\label{expg}Suppose $k\geq2$ and $X_{1},\ldots,X_{k}$ are finite nonempty
disjoint sets. Set
\[
\mathcal{PG}\left(  X_{1},\ldots,X_{k}\right)  =\left(  X_{1}\times
\cdots\times X_{k},2^{X_{1}\times\cdots\times X_{k}},\mu^{k},\mathbb{P}\left(
X_{1},\ldots,X_{k}\right)  \right)  ,
\]
where
\[
\mathbb{P}\left(  X_{1},\ldots,X_{k}\right)  =\left\{  A_{1}\times\cdots\times
A_{k}:A_{i}\subset X_{i}\text{ for all }i\in\left[  k\right]  \right\}
\]
and $\left(  X_{1}\times\cdots\times X_{k},2^{X_{1}\times\cdots\times X_{k}%
},\mu^{k}\right)  $ is the measure triple defined in Example \ref{exmp}. Then
$\mathcal{PG}\left(  X_{1},\ldots,X_{k}\right)  $ is an SR-system.
\end{example}

\subsubsection{Extending $\varepsilon$-regularity}

The primary goal of introducing SR-systems is to extend the concept of
$\varepsilon$-regular pairs\ of Szemer\'{e}di \cite{Sze76} (see also
\cite{Bol98}). For every $A,V\in\mathbb{A}$ set
\[
d\left(  A,V\right)  =\frac{\mu\left(  A\cap V\right)  }{\mu\left(  V\right)
}%
\]
if $\mu\left(  V\right)  >0,$ and $d\left(  A,V\right)  =0$ if $\mu\left(
V\right)  =0$.

\begin{definition}
\label{unis}Let $0<\varepsilon<1,\ V\in\mathbb{S},$ and $\mu\left(  V\right)
>0$. We call a set $A\in\mathbb{A}$ $\varepsilon$\textbf{-regular} in $V$ if
\[
\left\vert d\left(  A,U\right)  -d\left(  A,V\right)  \right\vert <\varepsilon
\]
for every $U\in\mathbb{S}$ such that $U\subset V$ and $\mu\left(  U\right)
>\varepsilon\mu\left(  V\right)  $.
\end{definition}

Let us see what Definition \ref{unis} says about directed $k$-graphs.

Take the SR-system $\mathcal{G}^{k}\left(  n\right)  .$ Let $G\in
\mathcal{G}^{k}\left(  n\right)  $ be a labelled directed $k$-graph with
$V\left(  G\right)  =\left[  n\right]  ,$ and let $\left(  V_{1},\ldots
,V_{k}\right)  $ be an ordered $k$-tuple of disjoint nonempty subsets of
$\left[  n\right]  .$ Write $e\left(  V_{1},\ldots,V_{k}\right)  $ for the
number of edges $\left(  v_{1},\ldots,v_{k}\right)  \in E\left(  G\right)  $
such that $v_{i}\in V_{i}$ for $i=1,\ldots,k.$

Now if $G$ is $\varepsilon^{1/k}$-regular in $V_{1}\times\cdots\times V_{k},$
then, for every ordered $k$-tuple $\left(  U_{1},\ldots,U_{k}\right)  $ such
that $U_{i}\subset V_{i}$ and $\left\vert U_{i}\right\vert >\varepsilon
\left\vert V_{i}\right\vert $ for $i=1,\ldots,k,$ we obtain%
\[
\left\vert \frac{e\left(  V_{1},\ldots,V_{k}\right)  }{\left\vert
V_{1}\right\vert \cdots\left\vert V_{k}\right\vert }-\frac{e\left(
U_{1},\ldots,U_{k}\right)  }{\left\vert U_{1}\right\vert \cdots\left\vert
U_{k}\right\vert }\right\vert <\varepsilon.
\]

Note that for $k=2$ this condition is essentially equivalent to the definition
of an \textquotedblleft$\varepsilon$-regular pair\textquotedblright.

Finally let us define $\varepsilon$-regularity with respect to partitions.

\begin{definition}
\label{unip}Let $0<\varepsilon<1$ and $\mathcal{P}$ be a partition of
$X\mathbb{\ }$into sets belonging to $\mathbb{S}$. We call a set
$A\in\mathbb{A}$ $\varepsilon$-\textbf{regular} in $\mathcal{P}$ if
\[%
{\textstyle\sum}
\left\{  \mu\left(  P\right)  :P\in\mathcal{P},\text{ }A\text{ is not
}\varepsilon\text{-regular in }P\right\}  <\varepsilon.
\]

\end{definition}

\subsection{Partitions in measure triples}

Given a collection $\mathbb{X}$ of subsets of $X,$ we write $\Pi\left(
\mathbb{X}\right)  $ for the family of finite partitions of $X$ into sets
belonging to $\mathbb{X}.$ We shall be mainly interested in $\Pi\left(
\mathbb{S}\right)  .$

Let $\mathcal{P},$ $\mathcal{Q}$ be partitions of $X,$ and $A\subset X.$ We
say that $\mathcal{P}$ refines $A$ (in notation $\mathcal{P}\succ A$) if $A$
is a union of\ members of $\mathcal{P},$ and that $\mathcal{P}$ refines
$\mathcal{Q}$ (in notation $\mathcal{P}\succ\mathcal{Q}$) if $\mathcal{P}$
refines each $Q\in\mathcal{Q}.$ We write $\mathcal{P}\cap\mathcal{Q}$ for the
partition consisting of all nonempty intersections $P\cap Q,$ where
$P\in\mathcal{P}$ and $Q\in\mathcal{Q}.$

Define the partition $\mathcal{P}^{k}$ of $X^{k}$ as%
\[
\mathcal{P}^{k}=\left\{  P_{i_{1}}\times\cdots\times P_{i_{k}}:P_{i_{j}}%
\in\mathcal{P},\text{ for all }j\in\left[  k\right]  \right\}  .
\]

\subsubsection{Bounding families of partitions}

We say that a family of partitions $\Phi\subset\Pi\left(  \mathbb{S}\right)  $
\emph{bounds} $\Pi\left(  \mathbb{S}\right)  $ if for every $\mathcal{P}\in
\Pi\left(  \mathbb{S}\right)  ,$ there exists $\mathcal{Q}\in\Phi$ such that
$\mathcal{Q}\succ\mathcal{P}$ and $\left\vert \mathcal{Q}\right\vert
\leq\varphi\left(  \left\vert \mathcal{P}\right\vert \right)  ,$ where
$\varphi:%
\mathbb{N}
\rightarrow%
\mathbb{N}
$ is a fixed increasing function, the \emph{rate }of $\Phi$.\bigskip

Here is an example of a bounding family. Given an integer $k\geq2$, take the
SR-system $\left(  X^{k},a\left(  \mathbb{A}^{k}\right)  ,\mu^{k}%
,\mathbb{S}^{\left\langle k\right\rangle }\right)  $ and define a family of
partitions $\Phi^{k}\subset\Pi\left(  \mathbb{S}^{\left\langle k\right\rangle
}\right)  $ as
\[
\Phi^{k}=\left\{  \mathcal{F}^{k}:\mathcal{F}\in\Pi\left(  \mathbb{S}\right)
\right\}  .
\]

\begin{lemma}
\label{le2} The family $\Phi^{k}$ bounds $\Pi\left(  \mathbb{S}^{\left\langle
k\right\rangle }\right)  $.
\end{lemma}

\section{The main result}

We are ready now to state our main theorem whose proof is presented in
\ref{pmt}.

\begin{theorem}
\label{gmaint} Let $\left(  X,\mathbb{A},\mu,\mathbb{S}\right)  $ be an
SR-system and suppose $\Phi$ is a family of partitions bounding $\Pi\left(
\mathbb{S}\right)  .$ Given a finite collection of measurable sets
$\mathcal{L}\subset\mathbb{A},$ a partition $\mathcal{P}\in\Pi\left(
\mathbb{S}\right)  ,$ and $\epsilon>0$, there exists $q=q\left(
\epsilon,\left\vert \mathcal{L}\right\vert ,\left\vert \mathcal{P}\right\vert
\right)  $ and $\mathcal{Q}\in\Phi$ such that:

- $\mathcal{Q}\succ\mathcal{P};$

- every $A\in\mathcal{L}$ is $\epsilon$-regular in $\mathcal{Q}$;

- $\left\vert \mathcal{Q}\right\vert \leq q.$
\end{theorem}

Our next goal is to show that Theorem \ref{gmaint} implies various types of
regularity lemmas. We emphasize the three steps that are necessary for its application:

\emph{(i)} select a measure triple $\left(  X,\mathbb{A},\mu\right)  $;

\emph{(ii)} introduce $\varepsilon$-regularity by fixing a boundedly built
semi-ring $\mathbb{S\subset A}$;

\emph{(iii)} select a bounding family of partitions $\Phi\subset\Pi\left(
\mathbb{S}\right)  $ by demonstrating an upper bound on its rate
$\varphi\left(  p\right)  .$

We turn now to specific applications.

\section{\label{sul}Examples}

To obtain more familiar versions of the Regularity Lemma, we extend the
concept of \textquotedblleft$\varepsilon$-equitable
partitions\textquotedblright\ and investigate when such partitions form
bounding families.

\subsection{Equitable partitions}

Given $\varepsilon>0$ and a measure triple $\left(  X,\mathbb{A},\mu\right)
,$ a partition $\mathcal{P}=\left\{  P_{0},\ldots,P_{p}\right\}  \in\Pi\left(
\mathbb{A}\right)  $ is called $\varepsilon$\emph{-equitable}, if $\mu\left(
P_{0}\right)  \leq\varepsilon$ and $\mu\left(  P_{1}\right)  =\cdots
=\mu\left(  P_{p}\right)  \leq\varepsilon.$

Let $k\geq2,$ take the SR-system $\left(  X^{k},\mathcal{A}\left(
\mathbb{A}^{k}\right)  ,\mu^{k},\mathbb{S}^{\left\langle k\right\rangle
}\right)  ,$ and define a family of partitions $\Phi^{k}\left(  \varepsilon
\right)  \subset\Pi\left(  \mathbb{S}^{\left\langle k\right\rangle }\right)  $
as%
\begin{equation}
\Phi^{k}\left(  \varepsilon\right)  =\left\{  \mathcal{P}^{k}:\mathcal{P}%
\in\Pi\left(  \mathbb{A}\right)  \text{ and }\mathcal{P}\text{ is }%
\varepsilon\text{-equitable}\right\}  . \label{defFe}%
\end{equation}

It is possible to prove that under some mild conditions on $\left(
X,\mathbb{A},\mu\right)  $ the family $\Phi^{k}\left(  \varepsilon\right)  $
bounds $\Pi\left(  \mathbb{S}^{\left\langle k\right\rangle }\right)  .$ To
avoid technicalities, we prove this claim for the SR-system $\mathcal{G}%
^{k}\left(  n\right)  =\left(  \left[  n\right]  ^{k},2^{\left[  n\right]
^{k}},\mu^{k},\left(  2^{\left[  n\right]  }\right)  ^{\left\langle
k\right\rangle }\right)  $. Let $\Phi^{k}\left(  n,\varepsilon\right)  $ be
defined by (\ref{defFe}) for the SR-system $\mathcal{G}^{k}\left(  n\right)
.$

\begin{lemma}
\label{leref}Let $0<\varepsilon<1$ and $n>1/\varepsilon.$ The family $\Phi
^{k}\left(  n,\varepsilon\right)  $ bounds $\Pi\left(  \left(  2^{\left[
n\right]  }\right)  ^{\left\langle k\right\rangle }\right)  $ with rate%
\[
\varphi\left(  p\right)  =\left(  \left\lceil 2/\varepsilon\right\rceil
+1\right)  ^{k}2^{pk^{2}}.
\]

\end{lemma}

Likewise, let $\Phi^{k}\left(  \left[  0,1\right]  ,\varepsilon\right)  $ be
defined by \ref{defFe} for the SR-system $\mathcal{B}^{k}=\left(  \left[
0,1\right]  ^{k},\mathbb{B}^{k},\lambda^{k},\mathbb{B}^{\left\langle
k\right\rangle }\right)  .$

\begin{lemma}
\label{leref1}Let $0<\varepsilon<1.$ The family $\Phi^{k}\left(
\varepsilon\right)  $ bounds $\Pi\left(  \mathbb{B}^{\left\langle
k\right\rangle }\right)  $ with rate%
\[
\varphi\left(  p\right)  =\left(  \left\lceil 1/\varepsilon\right\rceil
+1\right)  ^{k}2^{pk^{2}}.
\]

$\hfill\square$
\end{lemma}

\subsubsection{Regularity lemmas for $k$-graphs}

We first state a Regularity Lemma for directed $k$-graphs. As noted above we
represent directed $k$-graphs as subsets of $2^{\left[  n\right]  ^{k}}$ and
define regularity in terms of the SR-system $\mathcal{G}^{k}\left(  n\right)
=\left(  \left[  n\right]  ^{k},2^{\left[  n\right]  ^{k}},\mu^{k},\left(
2^{\left[  n\right]  }\right)  ^{\left\langle k\right\rangle }\right)  $.

\begin{theorem}
\label{gmaint1} For all $0<\varepsilon<1$ and positive integers $k,l,$ there
exist $n_{0}\left(  k,\varepsilon\right)  $ and $q\left(  k,l,\varepsilon
\right)  $ such that if $n>n_{0}\left(  k,\varepsilon\right)  $ and
$\mathcal{L}$ is a collection of $l$ subsets directed $k$-graphs on the vertex
set $\left[  n\right]  ,$ then there exists a partition $\mathcal{Q}=\left\{
Q_{0},\ldots,Q_{q}\right\}  $ of $\left[  n\right]  $ satisfying

\emph{(i)} $q\leq q\left(  k,l,\varepsilon\right)  $;

\emph{(ii)} $\left\vert Q_{0}\right\vert <\varepsilon n$, $\left\vert
Q_{1}\right\vert =\cdots=\left\vert Q_{q}\right\vert <\varepsilon n;$

\emph{(iii)} Every graph $G\in\mathcal{L}$ is $\varepsilon$-regular in at
least $\left(  1-\varepsilon\right)  q^{k}$ sets $Q_{i_{1}}\times\cdots\times
Q_{i_{k}}$, where $\left(  i_{1},\ldots,i_{k}\right)  $ is a $k$-tuple of
distinct elements of $\left[  q\right]  .$
\end{theorem}

\bigskip

As a consequence we obtain a Regularity Lemma for undirected $k$-graphs. For
$k=2$ this is the result of Szemer\'{e}di, for $k>2$ this is the result of
Chung \cite{Chu1}. Recall that undirected $k$-graphs are subsets
$G\subset2^{\left[  n\right]  ^{k}}$ such that if $\left(  v_{1},\ldots
,v_{k}\right)  \in G,$ then $\left\{  v_{1},\ldots,v_{k}\right\}  $ is a set
of size $k$ and $G$ contains each permutation of $\left(  v_{1},\ldots
,v_{k}\right)  $.

\begin{theorem}
\label{gmaint1.1} For all $0<\varepsilon<1$ and positive integers $k,l,$ there
exist $n_{0}\left(  k,\varepsilon\right)  $ and $q\left(  k,l,\varepsilon
\right)  $ such that if $n>n_{0}\left(  k,\varepsilon\right)  $ and
$\mathcal{L}$ is a collection of $l$ undirected $k$-graphs on the vertex set
$\left[  n\right]  ,$ then there exists a partition $\mathcal{Q}=\left\{
Q_{0},\ldots,Q_{q}\right\}  $ of $\left[  n\right]  $ satisfying:

\emph{i)} $q\leq q\left(  k,l,\varepsilon\right)  $;

\emph{ii)} $\left\vert Q_{0}\right\vert <\varepsilon n$, $\left\vert
Q_{1}\right\vert =\cdots=\left\vert Q_{q}\right\vert <\varepsilon n;$

\emph{iii)} For every graph $G\in\mathcal{L},$ there exist at least $\left(
1-\varepsilon\right)  \binom{q}{k}$ sets $\left\{  i_{1},\ldots,i_{k}\right\}
$ of distinct elements of $\left[  q\right]  $ such that $G$ is $\varepsilon
$-regular in $Q_{j_{1}}\times\cdots\times Q_{j_{k}}$ for every permutation of
$\left(  j_{1},\ldots,j_{k}\right)  $ of $\left\{  i_{1},\ldots,i_{k}\right\}
.$
\end{theorem}

\subsubsection{A regularity lemma for $k$-partite $k$-graphs}

Considering the SR-system $\mathcal{PG}\left(  X_{1},\ldots,X_{k}\right)  $
from Example \ref{expg} we obtain a regularity lemma for $k$-partite $k$-graphs.

\begin{theorem}
\label{gmaint1.2} Let $X_{1},\ldots,X_{k}$ be disjoint sets with $\left\vert
X_{1}\right\vert =\cdots=\left\vert X_{k}\right\vert =n.$ For all
$0<\varepsilon<1$ and positive integers $k,l,$ there exist $n_{0}\left(
k,\varepsilon\right)  $ and $q\left(  k,l,\varepsilon\right)  $ such that if
$n>n_{0}\left(  k,\varepsilon\right)  $ and $\mathcal{L}$ is a collection of
$l$ undirected $k$-partite $k$-graphs with vertex classes $X_{1},\ldots
,X_{k},$ then for each $i\in\left[  k\right]  ,$ there exist a partition
$\mathcal{Q}_{i}=\left\{  Q_{i0},\ldots,Q_{iq}\right\}  $ of $X_{i},$ satisfying:

\emph{i)} $q\leq q\left(  k,l,\varepsilon\right)  $;

\emph{ii)} $\left\vert Q_{i,0}\right\vert <\varepsilon n$, $\left\vert
Q_{i,1}\right\vert =\cdots=\left\vert Q_{i,q}\right\vert <\varepsilon n;$

\emph{iii)} For every graph $G\in\mathcal{L},$ there exist at least $\left(
1-\varepsilon\right)  q^{k}$ vectors $\left(  i_{1},\ldots,i_{k}\right)
\in\left[  q\right]  ^{k}$ such that $G$ is $\varepsilon$-regular in
$Q_{1,i_{1}}\times\cdots\times Q_{k,i_{k}}.$\bigskip
\end{theorem}

\subsubsection{A regularity lemma for measurable subsets of the unit cube}

Now define regularity according to the SR-system $\mathcal{B}^{k}=\left(
\left[  0,1\right]  ^{k},\mathbb{B}^{k},\lambda^{k},\mathbb{B}^{\left\langle
k\right\rangle }\right)  $. We obtain the following result.

\begin{theorem}
\label{gmaint2} For all $0<\varepsilon<1$ and positive integers $k,l,$ there
exists $q\left(  k,l,\varepsilon\right)  $ such that if $\mathcal{L}$ is a
collection of $l$ measurable subsets of the cube $\left[  0,1\right]  ^{k},$
then there exists a partition $\mathcal{Q}=\left\{  Q_{0},\ldots
,Q_{q}\right\}  $ of $\left[  0,1\right]  $ into measurable sets satisfying:

\emph{i)} $q\leq q\left(  k,l,\varepsilon\right)  $;

\emph{ii)} $\mu\left(  Q_{0}\right)  <\varepsilon$, $\mu\left(  Q_{1}\right)
=\cdots=\mu\left(  Q_{q}\right)  <\varepsilon;$

\emph{iii)} Every set $L\in\mathcal{L}$ is $\varepsilon$-regular in at least
$\left(  1-\varepsilon\right)  q^{k}$ sets $Q_{i_{1}}\times\cdots\times
Q_{i_{k}}$, where $\left(  i_{1},\ldots,i_{k}\right)  $ is a $k$-tuple of
distinct elements of $\left[  q\right]  .$
\end{theorem}

\bigskip

Finally let us define regularity according to the SR-system $\mathcal{B}%
^{k}=\left(  \left[  0,1\right]  ^{k},\mathbb{B}^{k},\lambda^{k}%
,\mathbb{I}^{\left\langle k\right\rangle }\right)  .$ We obtain a result which
we believe is specific to our approach.

\begin{theorem}
For all $0<\varepsilon<1$ and positive integers $k,l,$ there exists $q\left(
k,l,\varepsilon\right)  $ such that if $\mathcal{L}$ is a collection of $l$
measurable subsets of the cube $\left[  0,1\right]  ^{k}$ then there exists a
partition $\mathcal{Q}=\left\{  Q_{0},\ldots,Q_{q}\right\}  $ of $\left[
0,1\right]  $ satisfying:

\emph{i)} $q\leq q\left(  k,l,\varepsilon\right)  $;

\emph{ii)} $\mu\left(  Q_{0}\right)  <\varepsilon$, and the sets $Q_{1}%
,\ldots,$ $Q_{q}$ are intervals of equal length $<\varepsilon;$

\emph{iii)} Every set $L\in\mathcal{L}$ is $\varepsilon$-regular in at least
$\left(  1-\varepsilon\right)  q^{k}$ bricks $Q_{i_{1}}\times\cdots\times
Q_{i_{k}}$, where $\left(  i_{1},\ldots,i_{k}\right)  $ is a $k$-tuple of
distinct elements of $\left[  q\right]  .$
\end{theorem}

\section{\label{sprf}Proofs}

\subsection{\label{pmt}Proof of Theorem \ref{gmaint}}

Our proof is an adaptation of the original proof of Szemer\'{e}di \cite{Sze76}
(see also \cite{Bol98}). The following basic lemma is known as the
\textquotedblleft defect form of the Cauchy-Schwarz
inequality\textquotedblright; for a proof see \cite{Bol98}.

\begin{lemma}
\label{clem1} Let $x_{i}$ and $c_{i}$ be positive numbers for $i=1,\ldots,n.$
Then
\[%
{\textstyle\sum\limits_{i=1}^{n}}
c_{i}%
{\textstyle\sum\limits_{i=1}^{n}}
c_{i}x_{i}^{2}\geq\left(
{\textstyle\sum\limits_{i=1}^{n}}
c_{i}x_{i}\right)  ^{2}.
\]
If\ $J$ is a proper subset of $\left[  n\right]  $ and $\gamma>0$ is such
that
\[%
{\textstyle\sum\limits_{i=1}^{n}}
c_{i}%
{\textstyle\sum\limits_{i\in J}}
c_{i}x_{i}\geq%
{\textstyle\sum\limits_{i=1}^{n}}
c_{i}x_{i}%
{\textstyle\sum\limits_{i\in J}}
c_{i}+\gamma,
\]
then
\[%
{\textstyle\sum\limits_{i=1}^{n}}
c_{i}%
{\textstyle\sum\limits_{i=1}^{n}}
c_{i}x_{i}^{2}\geq\left(
{\textstyle\sum\limits_{i=1}^{n}}
c_{i}x_{i}\right)  ^{2}+\gamma^{2}/\left(
{\textstyle\sum\limits_{i\in J}}
c_{i}%
{\textstyle\sum\limits_{i\in\left[  n\right]  \backslash J}}
c_{i}\right)  .
\]

$\hfill\square$
\end{lemma}

Let $\mathcal{P}=\left\{  P_{1},\ldots,P_{p}\right\}  \in\Pi\left(
\mathbb{S}\right)  ,$ and $A\in\mathbb{A}$. Define the \emph{index} of
$\mathcal{P}$ with respect to $A$ as
\[
ind_{A}\mathcal{P}=%
{\textstyle\sum\limits_{P_{i}\in\mathcal{P}}}
\mu\left(  P_{i}\right)  d^{2}\left(  A\cap P_{i}\right)  .
\]
Note that for every $A\in\mathbb{A},$%
\begin{equation}
ind_{A}\mathcal{P}=%
{\textstyle\sum\limits_{P_{i}\in\mathcal{P}}}
\mu\left(  P_{i}\right)  d^{2}\left(  A\cap P_{i}\right)  \leq%
{\textstyle\sum\limits_{P_{i}\in\mathcal{P}\mathbf{,}\text{ }\mu\left(
P_{i}\right)  >0}}
\frac{\mu\left(  A\cap P_{i}\right)  \mu\left(  P_{i}\right)  }{\mu\left(
P_{i}\right)  }=\mu\left(  A\right)  \leq1. \label{mxind}%
\end{equation}

\begin{lemma}
\label{clem2} If $\mathcal{P},\mathcal{Q}\in\Pi\left(  \mathbb{S}\right)  ,$
$A\in\mathbb{A},$ and $\mathcal{Q}\succ\mathcal{P}$ then $ind_{A}%
\mathcal{Q}\geq ind_{A}\mathcal{P}$.
\end{lemma}

\begin{proof}
For simplicity we shall assume that $\mathcal{P}$ and $\mathcal{Q}$ consist
only of sets of positive measure. Fix $P\in\mathcal{P}$ and for every
$Q_{i}\subset P$, set
\[
c_{i}=\mu\left(  Q_{i}\right)  \text{ \ \ and \ }x_{i}=d\left(  A\cap
Q_{i}\right)  .
\]
Note that
\[%
{\textstyle\sum\limits_{Q_{i}\subset P}}
c_{i}=%
{\textstyle\sum\limits_{Q_{i}\subset P}}
\mu\left(  Q_{i}\right)  =\mu\left(  P\right)  \text{ \ \ and \ \ }\sum
_{Q_{i}\subset P}c_{i}x_{i}=\mu\left(  A\cap P\right)  .
\]
The Cauchy-Schwarz inequality (the first part of Lemma \ref{clem1}) implies
that%
\[%
{\textstyle\sum\limits_{Q_{i}\subset P}}
\mu\left(  Q_{i}\right)  d^{2}\left(  A\cap Q_{i}\right)  =%
{\textstyle\sum\limits_{Q_{i}\subset P}}
c_{i}x_{i}^{2}\geq\frac{1}{\mu\left(  P\right)  }\left(
{\textstyle\sum\limits_{Q_{i}\subset P}}
c_{i}x_{i}\right)  ^{2}=\frac{\mu^{2}\left(  A\cap P\right)  }{\mu\left(
P\right)  }.
\]
Summing over all sets $P\in\mathcal{P}$, the desired inequality follows.
\end{proof}

The next lemma supports the proof of Lemma \ref{clem4}.

\begin{lemma}
\label{clem3} Suppose $A,S,T\in\mathbb{A},$ $T\subset S$ and $\mu\left(
T\right)  >0.$ If
\begin{equation}
\left\vert d\left(  A\cap T\right)  -d\left(  A\cap S\right)  \right\vert
\geq\epsilon\label{ineq3}%
\end{equation}
then every partition $\mathcal{U}=\left\{  U_{1},\ldots,U_{p}\right\}  \in
\Pi\left(  \mathbb{A}\right)  $ such that $\mathcal{U}\succ S$ and
$\mathcal{U}\succ T,$ satisfies
\[%
{\textstyle\sum\limits_{U_{i}\subset S\text{ }}}
\mu\left(  U_{i}\right)  d^{2}\left(  A,U_{i}\right)  \geq\mu\left(  S\right)
d^{2}\left(  A,S\right)  +\epsilon^{2}\mu\left(  T\right)  .
\]

\end{lemma}

\begin{proof}
Let the partition $\mathcal{U}=\left\{  U_{1},\ldots,U_{p}\right\}  \in
\Pi\left(  \mathbb{A}\right)  $ be such that $\mathcal{U}\succ S$ and
$\mathcal{U}\succ T.$ For every $U_{i}\subset S,$ set
\[
c_{i}=\mu\left(  U_{i}\right)  ,\text{ \ \ }x_{i}=d\left(  A,U_{i}\right)  ,
\]
and observe that%
\[%
{\textstyle\sum\limits_{U_{i}\subset S}}
c_{i}=%
{\textstyle\sum\limits_{U_{i}\subset S}}
\mu\left(  U_{i}\right)  =\mu\left(  S\right)  \text{ \ \ and \ \ }%
{\textstyle\sum\limits_{U_{i}\subset S}}
c_{i}x_{i}=%
{\textstyle\sum\limits_{U_{i}\subset S}}
\mu\left(  A\cap U_{i}\right)  =\mu\left(  A\cap S\right)  .
\]
Similarly, we have
\[%
{\textstyle\sum\limits_{U_{i}\subset T}}
c_{i}=\mu\left(  T\right)  \text{ \ \ and \ \ }%
{\textstyle\sum\limits_{U_{i}\subset T}}
c_{i}x_{i}=%
{\textstyle\sum\limits_{U_{i}\subset T}}
\mu\left(  A\cap U_{i}\right)  =\mu\left(  A\cap T\right)  .
\]

Inequality (\ref{ineq3}) implies that either
\begin{equation}
d\left(  A,T\right)  >d\left(  A,S\right)  +\epsilon\text{ } \label{ineq4}%
\end{equation}
or
\[
d\left(  A,S\right)  >d\left(  A,T\right)  +\epsilon.
\]

Assume that (\ref{ineq4}) holds; the argument in the other case is identical.
Hence, $\mu\left(  T\right)  \neq\mu\left(  S\right)  ,$ so $T\subset S$
implies that $\mu\left(  T\right)  <\mu\left(  S\right)  .$ Furthermore,%
\begin{align*}%
{\textstyle\sum\limits_{U_{i}\subset T}}
c_{i}x_{i}  &  =d\left(  A,T\right)  \mu\left(  T\right)  >\left(  d\left(
A,S\right)  +\epsilon\right)  \mu\left(  T\right)  =\left(  d\left(
A,S\right)  +\epsilon\right)  \sum_{U_{i}\subset T}c_{i}\\
&  =\frac{%
{\textstyle\sum\nolimits_{U_{i}\subset S}}
c_{i}x_{i}}{\sum_{U_{i}\subset S}c_{i}}%
{\textstyle\sum\limits_{U_{i}\subset T}}
c_{i}+\epsilon%
{\textstyle\sum\limits_{U_{i}\subset T}}
c_{i}.
\end{align*}
By the definition of $c_{i}$ and $x_{i},$ we have
\[%
{\textstyle\sum\limits_{U_{i}\subset S}}
\mu\left(  U_{i}\right)  d^{2}\left(  A,U_{i}\right)  =%
{\textstyle\sum\limits_{U_{i}\subset S}}
c_{i}x_{i}^{2}.
\]
Therefore, setting
\[
\lambda=\epsilon%
{\textstyle\sum\limits_{U_{i}\subset S}}
c_{i}%
{\textstyle\sum\limits_{U_{i}\subset T}}
c_{i}=\epsilon\mu\left(  T\right)  \mu\left(  S\right)  ,
\]
and applying the second part of Lemma \ref{clem1}, we find that
\begin{align*}
\mu\left(  S\right)
{\textstyle\sum\limits_{U_{i}\subset S}}
\mu\left(  U_{i}\right)  d^{2}\left(  A,U_{i}\right)   &  \geq\left(
{\textstyle\sum\limits_{U_{i}\subset S}}
c_{i}x_{i}\right)  ^{2}+\lambda^{2}/\left(
{\textstyle\sum\limits_{U_{i}\subset T}}
c_{i}%
{\textstyle\sum\limits_{U_{i}\subset S\backslash T}}
c_{i}\right) \\
&  =\mu^{2}\left(  A\cap S\right)  +\epsilon^{2}\frac{\mu^{2}\left(  S\right)
\mu\left(  T\right)  }{\mu\left(  S\right)  -\mu\left(  T\right)  }.
\end{align*}
Hence,
\[%
{\textstyle\sum\limits_{U_{i}\subset S}}
{\textstyle\sum\limits_{U_{i}\subset S}}
\mu\left(  U_{i}\right)  d^{2}\left(  A,U_{i}\right)  \geq\mu\left(  S\right)
d^{2}\left(  A,S\right)  +\epsilon^{2}\frac{\mu\left(  T\right)  }{\mu\left(
S\right)  \left(  \mu\left(  S\right)  -\mu\left(  T\right)  \right)  }%
>\mu\left(  S\right)  d^{2}\left(  A,S\right)  +\epsilon^{2}\mu\left(
T\right)
\]
and this is exactly the desired inequality.
\end{proof}

The following lemma gives a condition for an absolute increase of
$ind_{A}\mathcal{P}$ resulting from refining.

\begin{lemma}
\label{clem4} Let $0<\epsilon<1$ and $\mathbb{S}$ be $r$-built. If
$\mathcal{P}\in\Pi\left(  \mathbb{S}\right)  $ and $A\in\mathbb{A}$ is not
$\epsilon$-regular in $\mathcal{P}$ then there exists $\mathcal{Q}\in
\Pi\left(  \mathbb{S}\right)  $ satisfying $\mathcal{Q}\succ\mathcal{P},$
$\left\vert \mathcal{Q}\right\vert \leq\left(  r+1\right)  \left\vert
\mathcal{P}\right\vert ,$ and
\begin{equation}
ind_{A}\mathcal{Q}\geq ind_{A}\mathcal{P}+\epsilon^{4}. \label{ineq1}%
\end{equation}

\end{lemma}

\begin{proof}
Let $\mathcal{P}=\left\{  P_{1},\ldots,P_{p}\right\}  $ and $\mathcal{N}$ be
the set of all $P_{i}$ for which $A$ is not $\epsilon$-regular in $P_{i}$.
Since $A$ is not $\epsilon$-regular in $\mathcal{P}$, by definition, we have%
\[%
{\textstyle\sum\limits_{P_{i}\in\mathcal{N}}}
\mu\left(  P_{i}\right)  \geq\epsilon.
\]

For every $P_{i}\in\mathcal{N},$ since $A$ is not $\epsilon$-regular in
$P_{i},$ there is a set $T_{i}\subset P_{i}$ such that $T_{i}\in\mathbb{S}$,
$\mu\left(  T_{i}\right)  >\epsilon\mu\left(  P_{i}\right)  ,$ and
\[
\left\vert d\left(  A,P_{i}\right)  -d\left(  A,T_{i}\right)  \right\vert
\geq\epsilon.
\]

Since $\mathbb{S}$ is $r$-built, for every $P_{i}\in\mathcal{N},$ there is a
partition of $P_{i}\backslash T_{i}$ into $r$ disjoint sets $A_{i1},\ldots
A_{is}\in\mathbb{S}$; hence $\left\{  A_{i1},\ldots A_{is},T_{i}\right\}  $ is
a partition of $P_{i}$ into at most $\left(  r+1\right)  $ sets belonging to
$\mathbb{S}.$ Let $\mathcal{Q}$ be the collection of all sets $A_{i1},\ldots
A_{is},T_{i},$ where $P_{i}\in\mathcal{N},$ together with all sets $P_{j}%
\in\mathcal{P}\backslash\mathcal{N}.$ Clearly $\mathcal{Q}\in\Pi\left(
\mathbb{S}\right)  ;$ also, $\mathcal{Q}\succ T_{i}$ and $\mathcal{Q}\succ
P_{i}$ for every $P_{i}\in\mathcal{N},$ and
\[
\left\vert \mathcal{Q}\right\vert \leq\left(  r+1\right)  \left\vert
\mathcal{N}\right\vert +\left\vert \mathcal{P}\right\vert -\left\vert
\mathcal{N}\right\vert \leq\left(  r+1\right)  \left\vert \mathcal{P}%
\right\vert .
\]

Thus, to finish the proof, we have to prove (\ref{ineq1}). Let $\mathcal{Q}%
=\left\{  Q_{1},\ldots,Q_{q}\right\}  .$ For every $P_{k}\in\mathcal{N},$
Lemma \ref{clem3} implies that%
\begin{equation}%
{\textstyle\sum\limits_{Q_{i}\subset P_{k}}}
\mu\left(  Q_{i}\right)  d^{2}\left(  A,Q_{i}\right)  \geq\mu\left(
P_{k}\right)  d^{2}\left(  A,P_{k}\right)  +\epsilon^{2}\mu\left(
T_{k}\right)  \geq\mu\left(  P_{k}\right)  d^{2}\left(  A,P_{k}\right)
+\epsilon^{3}\mu\left(  P_{k}\right)  . \label{eq1}%
\end{equation}

For any $P_{k}\in\mathcal{P},$ Lemma \ref{clem2} implies that
\begin{equation}%
{\textstyle\sum\limits_{Q_{i}\subset P_{k}}}
\mu\left(  Q_{i}\right)  d^{2}\left(  A,Q_{i}\right)  \geq\mu\left(
P_{k}\right)  d^{2}\left(  A,P_{k}\right)  . \label{eq2}%
\end{equation}
Now, by (\ref{eq1}) and (\ref{eq2}), we obtain
\begin{align*}
ind_{A}\mathcal{Q}  &  =%
{\textstyle\sum\limits_{Q_{i}\subset\mathcal{Q}}}
\mu\left(  Q_{i}\right)  d^{2}\left(  A,Q_{i}\right)  \geq\sum_{P_{i}%
\in\mathcal{P}}\mu\left(  P_{k}\right)  d^{2}\left(  A,P_{i}\right)
+\epsilon^{3}%
{\textstyle\sum\limits_{P_{i}\in\mathcal{N}}}
\mu\left(  P_{k}\right) \\
&  \geq%
{\textstyle\sum\limits_{P_{i}\in\mathcal{P}}}
\mu\left(  P_{k}\right)  d^{2}\left(  A,P_{i}\right)  +\epsilon^{4}%
=ind_{A}\mathcal{P}+\epsilon^{4},
\end{align*}
completing the proof.
\end{proof}

\begin{proof}
[Proof of Theorem \ref{gmaint}]Suppose $\mathbb{S}$ is $r$-built, $\Phi$
bounds $\Pi\left(  \mathbb{S}\right)  $ with rate $\varphi\left(
\cdot\right)  $ and let $\left\vert \mathcal{P}\right\vert =p.$ Define a
function $\psi:%
\mathbb{N}
\rightarrow%
\mathbb{N}
$ by
\begin{align}
\psi\left(  1,p\right)   &  =p;\label{defpsi}\\
\psi\left(  s+1,p\right)   &  =\left(  r+1\right)  \varphi\left(  \psi\left(
s,p\right)  \right)  ,\text{\ \ for\ every }s>1.\nonumber
\end{align}
We shall show that the partition $\mathcal{Q}\in\Phi$ may be selected so that
$\left\vert \mathcal{Q}\right\vert \leq\psi\left(  l\left\lfloor \epsilon
^{-4}\right\rfloor ,p\right)  .$

Select first a partition $\mathcal{P}_{0}\in\Phi$ such that $\mathcal{P}%
_{0}\succ\mathcal{P}$ and $\left\vert \mathcal{P}_{0}\right\vert \leq
\varphi\left(  \left\vert \mathcal{P}\right\vert \right)  $. We build
recursively a sequence of partitions $\mathcal{P}_{1},\mathcal{P}_{2},\ldots$
satisfying%
\begin{align}
\mathcal{P}_{i+1}  &  \succ\mathcal{P}_{i},\label{bnd0}\\
\left\vert \mathcal{P}_{i+1}\right\vert  &  \leq\varphi\left(  \left(
r+1\right)  \left\vert \mathcal{P}_{i}\right\vert \right)  ,\label{bnd1}\\
\exists A_{i}  &  \in\mathcal{G}:ind_{A_{i}}\mathcal{P}_{i+1}\geq ind_{A_{i}%
}\mathcal{P}_{i}+\epsilon^{4} \label{bnd2}%
\end{align}
for every $i=0,1,\ldots$ The sequence is built according the following rule:
If all $A\in\mathcal{L}$ are $\epsilon$-regular in $\mathcal{P}_{i},$ then we
stop. Otherwise there exists $A_{i}\in\mathcal{L}$ that is not $\epsilon
$-regular in $\mathcal{P}_{i}.$ Then, by Lemma \ref{clem4}, there is a
partition $\mathcal{P}_{i}^{\prime}\in\Pi\left(  \mathbb{S}\right)  $ such
that
\begin{align*}
\mathcal{P}_{i}^{\prime}  &  \succ\mathcal{P}_{i},\text{ }\\
\left\vert \mathcal{P}_{i}^{\prime}\right\vert  &  \leq\left(  r+1\right)
\left\vert \mathcal{P}_{i}\right\vert ,\text{ }\\
ind_{A_{i}}\mathcal{P}_{i}^{\prime}  &  \geq ind_{A_{i}}\mathcal{P}%
_{i}+\epsilon^{4}.
\end{align*}
Since $\Phi$ bounds $\mathbb{S}$ with rate $\varphi$, there is a partition
$\mathcal{P}_{i+1}\in\Phi$ such that $\mathcal{P}_{i+1}\succ\mathcal{P}%
_{i}^{\prime}$ and $\left\vert \mathcal{P}_{i+1}\right\vert \leq\varphi\left(
\left\vert \mathcal{P}_{i}^{\prime}\right\vert \right)  .$ Hence,
(\ref{bnd0}), (\ref{bnd1}), and (\ref{bnd2}) hold.

Set $k=\left\lfloor \epsilon^{-4}\right\rfloor .$ If the sequence
$\mathcal{P}_{0},\mathcal{P}_{1},\ldots$ has more than $lk$ terms then, by the
pigeonhole principle, there exist a set $A\in\mathcal{L}$ and a sequence
$\mathcal{P}_{i_{1}},\ldots,\mathcal{P}_{i_{k+1}}$, such that
\[
ind_{A}\mathcal{P}_{i_{j+1}}\geq ind_{A}\mathcal{P}_{i_{j}}+\epsilon^{4}%
\]
for every $j=1,\ldots,k.$ Hence, we find that
\[
ind_{A}\mathcal{P}_{i_{k+1}}\geq ind_{A}\mathcal{P}_{i_{1}}+k\epsilon
^{4}>k\epsilon^{4}\geq1,
\]
contradicting (\ref{mxind}). Therefore, all $A\in\mathcal{L}$ are $\epsilon
$-regular in some partition $\mathcal{Q}=\mathcal{P}_{i}$. By (\ref{bnd1}),
$\left\vert \mathcal{Q}\right\vert \leq\psi\left(  l\left\lfloor \epsilon
^{-4}\right\rfloor ,p\right)  $, completing the proof.
\end{proof}

\subsection{Proof of lemma \ref{leref}}

\begin{proof}
Select a partition $\mathcal{P}=\left\{  P_{1},\ldots,P_{p}\right\}  \in
\Pi\left(  \left(  2^{\left[  n\right]  }\right)  ^{\left\langle
k\right\rangle }\right)  ;$ for every $i\in\left[  p\right]  $ let
\[
P_{i}=R_{i1}\times\cdots\times R_{ik},\text{ }R_{ij}\subset\left[  n\right]
,\text{ for }j\in\left[  k\right]  .
\]
Let
\[
\mathcal{R}=\cap_{i\in\left[  p\right]  ,j\in\left[  k\right]  }\left\{
R_{ij},X\backslash R_{ij}\right\}  .
\]
and set $r=\left\vert \mathcal{R}\right\vert .$ Clearly, $r\leq2^{pk}.$ Our
first goal is to find an $\varepsilon$-equitable partition $\mathcal{Q}%
\succ\mathcal{R}$ with
\[
\left\vert \mathcal{Q}\right\vert \leq\left(  \frac{2}{\varepsilon}+1\right)
2^{pk}.
\]

Suppose first that $n\geq2r/\varepsilon.$ To construct the required
$\mathcal{Q},$ partition every $R\in\mathcal{R}$ into sets of size
$\left\lfloor \varepsilon n/r\right\rfloor $ and a smaller residual set. The
measure of each member of $\mathcal{Q}$ that is not residual is at most
$\left\lfloor \varepsilon n/r\right\rfloor /n\leq\varepsilon.$ The total
measure of all residual sets is less than%
\[
\frac{\left\lfloor \varepsilon n/r\right\rfloor }{n}r\leq\varepsilon,
\]
thus, $\mathcal{Q}$ is an $\varepsilon$-equitable partition refining
$\mathcal{R}$. Since
\[
\left\vert \mathcal{Q}\right\vert \leq\frac{n}{\left\lfloor \varepsilon
n/r\right\rfloor }+r\leq\frac{2n}{\varepsilon n/r}+r=\left(  \frac
{2}{\varepsilon}+1\right)  r\leq\left(  \frac{2}{\varepsilon}+1\right)
2^{pk},
\]
$\mathcal{Q}$ has the required properties.

Let now $n<2r/\varepsilon$ and $\mathcal{Q}$ be the partition of $\left[
n\right]  $ into $n$ sets of size $1$. Since $\varepsilon>1/n,$ the partition
$\mathcal{Q}$ is $\varepsilon$-equitable and refines $\mathcal{R}$. Since%
\[
\left\vert \mathcal{Q}\right\vert =n<\frac{2}{\varepsilon}r\leq\left(
\frac{2}{\varepsilon}+1\right)  2^{pk},
\]
$\mathcal{Q}$ has the required properties.

To complete the proof observe that $\mathcal{Q}^{k}\in\Phi^{k}\left(
\varepsilon\right)  ,$ $\mathcal{Q}^{k}\succ\mathcal{R}^{k}\succ\mathcal{P},$
and
\[
\left\vert \mathcal{Q}^{k}\right\vert \leq\left\vert \mathcal{Q}\right\vert
^{k}\leq\left(  \frac{2}{\varepsilon}+1\right)  ^{k}2^{pk^{2}}.
\]
\bigskip
\end{proof}

\textbf{Acknowledgement}

\end{document}